\documentclass[11pt,hyp,]{nyjm}
\usepackage{hyperref}
\hypersetup{nesting=true,debug=true,naturalnames=true}

\usepackage{graphicx,amssymb}
\recdate{August 18, 2013}
\newcommand{\figref}[1]{\hyperlink{#1}{\ref*{fig:#1}}}

\let\<\langle
\let\>\rangle
\usepackage[all,pdf]{xy}
\UseComputerModernTips
\newcommand{\doi}[1]{doi:\,\href{http://dx.doi.org/#1}{#1}}
\newcommand{\mailurl}[1]{\email{\href{mailto:#1}{#1}}}
\let\uml\"
\let\dmo\DeclareMathOperator

\dmo{\val}{val}

\usepackage{mathtools}
\usepackage{amssymb}

\usepackage{tikz}
\usetikzlibrary{topaths}

\title[Degree Theorem in Auter space]{A combinatorial proof of the Degree Theorem in Auter space} 

\author[Robert McEwen]{\hyperlink{mcewen}{Robert McEwen}}  
\address{\hypertarget{mcewen}Ruckersville, VA 22968} 
\mailurl{mcewen.rob@gmail.com}

\author[Matthew C. B. Zaremsky]{\hyperlink{zaremsky}{Matthew C. B. Zaremsky}}
\address{\hypertarget{zaremsky}Department of Mathematical Sciences, Binghamton University, Binghamton, NY 13902}
\mailurl{zaremsky@math.binghamton.edu}

\thanks{The second named author gratefully acknowledges support from the SFB~$701$ of the DFG}



\newtheorem{theorem}{Theorem}[section]
\newtheorem{lemma}[theorem]{Lemma}

\newtheorem*{nonumtheorem}{Theorem}
\newtheorem{corollary}[theorem]{Corollary}
\newtheorem{proposition}[theorem]{Proposition}
\newtheorem{observation}[theorem]{Observation}

\theoremstyle{definition}
\newtheorem{definition}[theorem]{Definition}

\newtheorem{claim}{Claim}

\newcommand{\edge}{\varepsilon}

\newcommand{\defeq}{\mathbin{\vcentcolon =}}
\newcommand{\eqdef}{\mathbin{= \vcentcolon}}

\DeclareMathOperator{\Aut}{Aut}

\DeclareMathOperator{\lk}{lk}
\DeclareMathOperator{\st}{st}
\DeclareMathOperator{\dlk}{{\lk}{\downarrow}}
\DeclareMathOperator{\dst}{{\st}{\downarrow}}
\DeclareMathOperator{\BU}{BU}
\DeclareMathOperator{\SBU}{SBU}

\begin{document} 
 
\begin{abstract}  
 We use discrete Morse theory to give a new proof of Hatcher and Vogtmann's \emph{Degree Theorem} in Auter space $A_n$. There is a filtration of $A_n$ into subspaces $A_{n,k}$ using the \emph{degree} of a graph, and the Degree Theorem says that each $A_{n,k}$ is $(k-1)$-connected. This result is useful, for example to calculate stability bounds for the homology of $\Aut(F_n)$. The standard proof of the Degree Theorem is global in nature. Here we give a proof that only uses local considerations, and lends itself more readily to generalization.
\end{abstract} 
\maketitle
\tableofcontents

\section{Introduction}\label{sec:introsec}

In this note we provide an alternate proof of Hatcher and Vogtmann's \emph{Degree Theorem} in Auter space \cite{hatcher98}, using discrete Morse theory. The advantage of our proof is that it relies only on local data, and also lends itself more readily to certain generalizations. \emph{Auter space} $A_n$ is the space of rank-$n$ basepointed marked metric graphs. In \cite{hatcher98}, a measurement called the \emph{degree} of a graph was used to filter $A_n$ into highly connected sublevel sets $A_{n,k}$, which were then used to produce stability bounds for the rational and integral homology of $\Aut(F_n)$. The key result was:

\begin{nonumtheorem}[Degree Theorem]\cite{hatcher98}
 $A_{n,k}$ is $(k-1)$-connected.
\end{nonumtheorem}

The proof of the Degree Theorem in \cite{hatcher98} is done by globally deforming disks in $A_n$ via an iterated process. Our proof here uses discrete Morse theory, as in \cite{bestvina97}, to reduce the problem to a purely local one. First we shift focus to the \emph{spine} of Auter space, which we denote $L_n$. This is a combinatorial model for $A_n$ that is a deformation retract. We construct a \emph{height function} $h$ on $L_n$ that reduces the problem to asking whether the \emph{descending links} with respect to $h$ are highly connected. This is advantageous for being a local rather than global problem, and also lends itself more readily to generalization. For example a similar method has been used in \cite{zar_psymmauts} to get stability results for the groups $\Sigma\!\Aut_n^m$ of \emph{partially symmetric} automorphisms.

In Section~\ref{sec:setupsec} we describe the spine of Auter space $L_n$, and define the notion of the \emph{degree} $d_0$ of a graph. We use the degree to filter $L_n$ into sublevel sets $L_{n,k}$, as in \cite{hatcher98}. We then define a height function $h$ on $L_n$ refining $d_0$, and consider the descending links of vertices in $L_n$ with respect to $h$. The descending link of a vertex decomposes as a join of two complexes, called the \emph{d-down-link} and \emph{d-up-link}. In Section~\ref{downsec} we analyze the connectivity of the d-down-link, and in Section~\ref{sec:upsec} we do the same for the d-up-link. The upshot of this is Corollary~\ref{cor:desc_lk_conn}, that the descending links are all highly connected. From this we quickly obtain that $L_{n,k}$, and hence $A_{n,k}$ is $(k-1)$-connected; see Theorem~\ref{thrm:deg_thm}.

\subsubsection*{Acknowledgments}  The authors would like to thank Kai-Uwe Bux, who helped with a preliminary version of this paper, and Allen Hatcher for his comments and suggestions. Parts of this paper are based on results in the first named author's Ph.D. thesis \cite{mcewen_thesis}, done at the University of Virginia.

\section{Auter space, degree, and a height function}\label{sec:setupsec}

We begin by describing the \emph{spine of Auter space} $L_n$ introduced in \cite{hatcher98}. Let $R_n$ be the rose with $n$ edges, i.e., the graph with a single vertex $p_0$ and $n$ edges. Here by a \emph{graph} we always mean a finite connected one-dimensional CW-complex, with the usual notions of vertices and edges. If $\Gamma$ is a rank $n$ graph with basepoint vertex $p$, a homotopy equivalence $\rho \colon R_n\to\Gamma$ taking $p_0$ to $p$ is called a \emph{marking} on $\Gamma$. Two markings are equivalent if there is a basepoint-preserving homotopy between them. We only consider graphs such that $p$ is at least bivalent and all other vertices are at least trivalent. The spine $L_n$ of Auter space is then the complex of marked basepointed rank $n$ graphs $(\Gamma,p,\rho)$, up to equivalence of markings.

To be more precise, $L_n$ is a simplicial complex with a vertex for every equivalence class of triples $(\Gamma,p,\rho)$. An $r$-simplex is given by a chain of \emph{forest collapses} $\Gamma_r\stackrel{d_r}{\to}\Gamma_{r-1}\stackrel{d_{r-1}}{\to}\cdots\stackrel{d_1}{\to}\Gamma_0$ and markings $\rho_i \colon R_n\to\Gamma_i$ with the following diagram commuting up to homotopy.

$$\xymatrix{
\Gamma_r \ar[r]^{d_r}
& \Gamma_{r-1} \ar[r]^{d_{r-1}}
& \cdots \ar[r]^{d_2} & \Gamma_1 \ar[r]^{d_1} & \Gamma_0 \\
& & R_n \ar[ull]^{\rho_r} \ar[ul]_{\rho_{r-1}} \ar[ur]^{\rho_1} \ar[urr]_{\rho_0}
}$$
Here a \emph{forest collapse} or \emph{blow-down} $d \colon \Gamma\to\Gamma'$ is a (basepoint-preserving) homotopy equivalence of graphs that is given by collapsing each component of a forest $F$ in $\Gamma$ to a point. We will write the resulting graph as $\Gamma/F$. The reverse of a blow-down is, naturally, called a \emph{blow-up}.

Let $\Gamma$ be a graph with rank $n$, basepoint $p$ and vertex set $V(\Gamma)$. The \emph{degree} of $\Gamma$ can be defined as \[\displaystyle d_0(\Gamma)\defeq \sum_{\mathclap{p\neq v\in V(\Gamma)}}(\val(v)-2)\] or equivalently as $d_0(\Gamma)=2n-\val(p)$ \cite[Section~3]{hatcher98}. Here $\val(v)$ is the valency of $v$, that is the number of half-edges incident to $v$. This is sometimes called the ``degree'' of the vertex, but we have reserved this word for the degree of a graph.

\begin{definition}[Filtration by degree]\label{def:Lnk}
 For $k\geq0$, let $L_{n,k}$ be the subcomplex of $L_n$ spanned by vertices represented by triples $(\Gamma,p,\rho)$ with $d_0(\Gamma)\leq k$.
\end{definition}

The Degree Theorem says that $A_{n,k}$ is $(k-1)$-connected, and this is equivalent to $L_{n,k}$ being $(k-1)$-connected \cite[Section~5.1]{hatcher98}, which is what we will prove.

We now define some other measurements on $\Gamma$. For $v\in V(\Gamma)$ let $d(p,v)$ denote the minimum length of an edge path in $\Gamma$ from $v$ to $p$, and call $d(p,v)$ the \emph{level} of $v$. Here we are treating each edge in the graph as having length $1$. Define $\Lambda_i(\Gamma)\defeq \{v\in V(\Gamma)\mid d(p,v)= i\}$, $n_i(\Gamma)\defeq -|\Lambda_i(\Gamma)|$ and
$$d_i(\Gamma)\defeq \sum_{\mathclap{v\in V(\Gamma)\setminus\Lambda_i(\Gamma)}}(\val(v)-2)$$
for $i\geq0$. Note that $\Lambda_0(\Gamma)=\{p\}$, $n_0(\Gamma)=-1$, and $d_0(\Gamma)$ agrees with the definition of degree, so this is not an abuse of notation. Finally, define
$$h(\Gamma)\defeq (d_0(\Gamma),n_1(\Gamma),d_1(\Gamma),n_2(\Gamma),d_2(\Gamma),\dots)$$
to be the \emph{height} of the graph $\Gamma$, considered with the lexicographic ordering. This height function is a refinement of the degree function. Extend the definition of $h$ to the vertices of $L_n$ via $h(\Gamma,p,\rho)=h(\Gamma)$. For brevity, in the future we will often just refer to vertices in $L_n$ as being graphs, rather than equivalence classes of triples $(\Gamma,p,\rho)$.

\begin{observation}\label{obs:sublev}
 $L_{n,k}$ is the sublevel set of $L_n$ defined by the inequality
 $$h(\Gamma)\le (k,1,0,0,\dots) \text{.}$$
\end{observation}

\begin{proof}
 If $h(\Gamma)\le (k,1,0,0,\dots)$ then $d_0(\Gamma)\le k$. Now suppose $d_0(\Gamma)\le k$. If $d_0(\Gamma)<k$ then $h(\Gamma)< (k,1,0,0,\dots)$. If $d_0(\Gamma)=k$ then since $n_1(\Gamma)\le 0$ we have $h(\Gamma)< (k,1,0,0,\dots)$.
\end{proof}

Any blow-down necessarily increases some $n_i$ (that is, decreases some $|\Lambda_i|$), and so adjacent vertices in $L_n$ have different heights. Hence $h$ is a ``true'' height function, in the sense of \cite{bestvina97}. This, together with Observation~\ref{obs:sublev}, means that the connectivity of $L_{n,k}$ can be deduced by inspecting the \emph{descending links} with respect to $h$ of vertices in $L_n\setminus L_{n,k}$. For a vertex $\Gamma$ in $L_n$, the \emph{descending star} $\dst(\Gamma)$ with respect to $h$ is the set of simplices in the star of $\Gamma$ whose vertices other than $\Gamma$ all have strictly lower height than $\Gamma$. The \emph{descending link} $\dlk(\Gamma)$ is the set of faces of simplices in $\dst(\Gamma)$ that do not themselves contain $\Gamma$.

There are two types of vertices in $\dlk(\Gamma)$: those obtained from $\Gamma$ by a descending blow-up, and those obtained by a descending blow-down. Here we say that a blow-up or blow-down is \emph{descending} if the resulting graph has a lower height than the starting graph. Call the full subcomplex of $\dlk(\Gamma)$ spanned by vertices of the first type the \emph{d-up-link}, and the subcomplex spanned by vertices of the second type the \emph{d-down-link}. Any vertex in the d-up-link is related to every vertex in the d-down-link by a blow-down, so $\dlk(\Gamma)$ is the simplicial join of the d-up- and d-down-links.

If blowing down the forest $F$ is a descending blow-down, we will call the forest itself \emph{descending}, and similarly a forest can be ascending. It will be a good idea to describe precisely which forests in a graph are ascending and descending. For a forest $F$ in $\Gamma$ define $D(F)\defeq \min\{i\mid F$ has a vertex in $\Lambda_i\}$ to be the \emph{level} of $F$. If there is an edge path in $F$ from a vertex in $\Lambda_{D(F)}$ to another, distinct vertex in $\Lambda_{D(F)}$, we say that $F$ \emph{connects vertices in $\Lambda_{D(F)}$}.

\begin{lemma}\label{lem:which_forests_desc}
 If $F$ connects vertices in $\Lambda_{D(F)}$ then $F$ is ascending. Otherwise $F$ is descending.
\end{lemma}

\begin{proof}
 Let $i\defeq D(F)$. Blowing down $F$ does not change any $n_j$ or $d_j$ for $j<i$. If $F$ connects vertices in $\Lambda_i$, then blowing down $F$ increases $n_i$, so $F$ is ascending. If $F$ does not connect any vertices in $\Lambda_i$, then blowing down $F$ will not change $n_i$, but since each non-basepoint vertex of $\Gamma$ is at least trivalent, $d_i$ will be smaller in $\Gamma/F$ than in $\Gamma$, and so $F$ is descending.
\end{proof}

As a corollary to the proof we obtain:

\begin{corollary}\label{which_blowups_desc}
 A blow-up at a vertex $v\in\Lambda_i$ is descending if and only if it decreases $n_i$, that is increases $|\Lambda_i|$.\qed
\end{corollary}

\begin{figure}[htb]\hypertarget{desc_blowup}\leavevmode
  \centering
  \includegraphics[scale=.4]{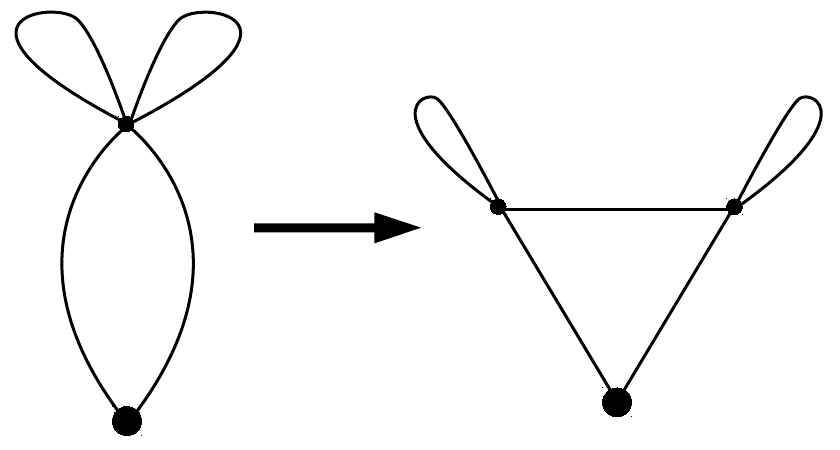}
  \caption{A descending blow-up.}
  \label{fig:desc_blowup}
\end{figure}

An example of a descending blow-up is given in Figure~\figref{desc_blowup}. Here $d_0$ stays constant $4$, and $n_1$ decreases from $-1$ to $-2$.

We close this section with some definitions regarding edges in graphs.

\begin{definition}\label{def:edges}
 Let $\edge$ be an edge in $\Gamma$, with vertices $v_1$ and $v_2$. We call $\edge$ \emph{horizontal} if $d(p,v_1)=d(p,v_2)$, and \emph{vertical} if $d(p,v_1)\neq d(p,v_2)$. Let $\edge$ be a vertical edge with vertices $v_1$ and $v_2$ such that $d(p,v_1)>d(p,v_2)$. We call $v_1$ the \emph{top} of $\edge$ and $v_2$ the \emph{bottom}. A half-edge can also have a top or a bottom (or neither, if it comes from a horizontal edge). We say that $\edge$ is \emph{decisive} if it is the unique vertical edge having $v_1$ as its top, that is if any minimal length edge path from $v_1$ to $p$ must begin with $\edge$.
 \end{definition}

\section{Connectivity of the d-down-link}\label{downsec}

In this section we analyze the d-down-link of $\Gamma$. In order for a certain induction to run, it will become necessary to consider (connected) graphs with vertices of valency $1$ and $2$. It turns out that $h$ does not ``work correctly'' on such graphs, for instance Lemma~\ref{lem:which_forests_desc} no longer holds. Therefore in this section we will use Lemma~\ref{lem:which_forests_desc} as a guide for which forests we want to consider.

Recall that we say $F$ connects vertices in $\Lambda_{D(F)}$ provided that there is an edge path in $F$ between distinct vertices of $\Lambda_{D(F)}$.

\begin{definition}
 Let $\Gamma$ be a connected graph with basepoint $p$, and with no restriction on the valency of vertices. Let $F$ be a subforest of $\Gamma$, with level $D(F)$. We will call $F$ \emph{bad} if it connects vertices in $\Lambda_{D(F)}$, and \emph{good} if it does not.
\end{definition}

Thanks to Lemma~\ref{lem:which_forests_desc}, if $\Gamma$ actually comes from $L_n$ then a forest in $\Gamma$ is good if and only if it is descending. Let $P(\Gamma)$ be the poset of good forests in $\Gamma$, ordered by inclusion, so if $\Gamma$ comes from $L_n$ then the geometric realization $|P(\Gamma)|$ of $P(\Gamma)$ is the d-down-link of $\Gamma$. Let $V$ be the number of vertices in $\Gamma$ and $E$ the number of edges. In what follows we will suppress the bars indicating geometric realization, so posets themselves will be said to have a homotopy type. Recall that an empty wedge of spheres is a single point.

\begin{proposition}[Homotopy type of the d-down-link]\label{prop:downprop}
 $P(\Gamma)$ is homotopy equivalent to a (possibly empty) wedge of spheres of dimension $V-2$.
\end{proposition}

\begin{proof}
 Our proof is similar to the proof of Proposition~2.2 in \cite{vogtmann90}. We induct on the number of edges $E$. We can assume that $\Gamma$ has no single-edge loops, since they do not affect $V$ or $P(\Gamma)$. We remark that already after this reduction the vertices may have arbitrary valency, so it is important that we are considering ``good'' forests instead of ``descending'' forests. Also, if $\Gamma$ has a separating edge $\edge$ then $P(\Gamma)$ is a cone with cone point $\edge$, so without loss of generality $\Gamma$ has no separating edges.
 
 The base case is $E=0$, for which $V=1$ and $P(\Gamma)=\emptyset=S^{V-2}$ as desired.

 Now suppose $E>0$. Choose an edge $\edge$ with endpoints $v_1,v_2$ maximizing the quantity $d(p,v_1)+d(p,v_2)$. In other words, $\edge$ is as far as possible from the basepoint; note that $D(\edge)$ is also maximized. Let $P_1(\Gamma)\subseteq P(\Gamma)$ be the poset of all good forests in $\Gamma$ except the forest just consisting of the edge $\edge$. Also let $P_0(\Gamma)\subseteq P_1(\Gamma)$ be the poset of good forests that do not contain $\edge$.
 
\medskip
 
 \begin{claim}\label{cl1} $P_1(\Gamma)\simeq P_0(\Gamma)$.
\end{claim}

\begin{proof}[Proof of Claim 1] For any $F\in P_1(\Gamma)$, $F-\edge$ is again a good forest by definition, so the poset map $g \colon P_1(\Gamma)\to P_1(\Gamma)$ given by $F\mapsto F-\edge$ is well defined. Here $F-\edge$ is just the forest obtained by removing $\edge$ from $F$. By construction, $g$ is the identity on its image $P_0(\Gamma)$, and $g(F)\leq F$ for all $F\in P_1(\Gamma)$, so $g$ induces a homotopy equivalence between $P_1(\Gamma)$ and $P_0(\Gamma)$ \cite[Section~1.3]{quillen78}. 
\end{proof}

 Now consider the graph $\Gamma-\edge$ obtained by removing $\edge$ from $\Gamma$. Since $\edge$ is not a separating edge, $\Gamma-\edge$ is connected.

 \begin{claim}\label{cl2} $P_0(\Gamma)\cong P(\Gamma-\edge)$.
\end{claim}

 \begin{proof}[Proof of Claim 2]  Consider the map $\iota \colon P(\Gamma-\edge)\to P_0(\Gamma)$ induced by $\Gamma-\edge\hookrightarrow\Gamma$. Since $D(\edge)$ is maximized and $\edge$ is not a separating edge, $\edge$ cannot be decisive, so adding $\edge$ to the graph does not change the levels $\Lambda_i$. In particular adding $\edge$ cannot affect whether a forest $F$ in $\Gamma-\edge$ is good or bad, so $\iota$ is an isomorphism. 
\end{proof}

 Since $\Gamma-\edge$ has $E-1$ edges and $V$ vertices, by induction $P(\Gamma-\edge)\simeq\bigvee S^{V-2}$. Then Claims \ref{cl1} and \ref{cl2} tell us that $P_1(\Gamma)\simeq\bigvee S^{V-2}$.

 With $P_1(\Gamma)$ in hand, we now ask about $P(\Gamma)$ itself. If $\edge$ is horizontal then it is bad, so $P_1(\Gamma)=P(\Gamma)$ and we are done. Assume instead that $\edge$ is vertical, hence good, which means $P(\Gamma)=P_1(\Gamma)\cup\st(\edge)$ with $P_1(\Gamma)\cap\st(\edge)=\lk(\edge)$, where link and star are taken in $P(\Gamma)$.
 
 Consider the graph $\Gamma/\edge$. This has $E-1$ edges and $V-1$ vertices, so by induction, $P(\Gamma/\edge)\simeq\bigvee S^{V-3}$. Hence it suffices now to prove the following:

 \begin{claim}\label{cl3} $\lk(\edge)\cong P(\Gamma/\edge)$.
\end{claim}

 \begin{proof}[Proof of Claim 3] First note that for a forest $F\neq\edge$ in $\Gamma$, $F$ is good if and only if $F/\edge$ is, where $F/\edge$ is the image of $F$ in $\Gamma/\edge$. Indeed, if $D(F)<D(\edge)$ then this is trivial; if $D(F)\geq D(\edge)$ then by our choice of $\edge$, $D(F)=D(\edge)$, and it is then evident that $F$ is good if and only if $F/\edge$ is. Now consider the map $c \colon \lk(\edge)\to P(\Gamma/\edge)$ sending $F$ to $F/\edge$. This is well-defined by the previous observation. We claim that $c$ is bijective. Let $\Phi\in P(\Gamma/\edge)$. There are precisely two forests in $\Gamma$ that map to $\Phi$ under blowing down $\edge$, one that contains $\edge$ and one that does not (this shows that $c$ is injective). Let $\Phi'$ be the one that does. If $\Phi$ was good then so is $\Phi'$, again by the previous observation, so $\Phi'\in\lk(\edge)$. Hence $c$ is an isomorphism. 
\end{proof}

This finishes the proof of the Proposition~\ref{prop:downprop}.
\end{proof}

It will also be convenient to establish one specific case when $P(\Gamma)$ is contractible.

\begin{lemma}\label{lem:downlemma2}
 If $\Gamma$ has a decisive edge then $P(\Gamma)$ is contractible.
\end{lemma}

\begin{proof}
 The proof is almost the same as the proof of the previous proposition. We again induct on $E$. If $E=0$ then $\Gamma$ does not have any edges, much less any decisive edges, and so the claim is vacuously true. Now assume $E>0$ and $\Gamma$ has a decisive edge $\eta$. If $\eta$ has maximum distance to the base point among edges in $\Gamma$ then it is separating and $P(\Gamma)$ is contractible with $\eta$ serving as a cone point. Otherwise, let $\edge\neq\eta$ be an edge in $\Gamma$ that has maximum distance to the basepoint, and define $P_1(\Gamma)$ and $P_0(\Gamma)$ as in the previous proof.

 By Claims \ref{cl1} and \ref{cl2} in the previous proof, $P_1(\Gamma)\simeq P_0(\Gamma)\cong P(\Gamma-\edge)$. This is contractible by induction since $\Gamma-\edge$ has fewer edges and still contains the decisive edge $\eta$. If $\edge$ is horizontal, $P(\Gamma)=P_1(\Gamma)$ and we are done, so assume $\edge$ is vertical. As in the previous proof, it then suffices to show that $\lk(\edge)$ has the appropriate homotopy type, i.e., is contractible. By Claim \ref{cl3} in the previous proof, $\lk(\edge)\simeq P(\Gamma/\edge)$. Let $\eta'$ be the image of $\eta$ in $\Gamma/\edge$. Since $\eta$ is decisive, $\edge$ and $\eta$ have different tops. Since $\edge$ is at maximal distance from $p$, $\eta'$ is a decisive edge in $\Gamma/\edge$. Hence $P(\Gamma/\edge)$ is contractible by induction, and we are done.
\end{proof}

\section{Connectivity of the d-up-link}\label{sec:upsec}

We now inspect the d-up-link. We first focus on one vertex at a time. Let $\BU(v)$ be the poset of all blow-ups at the vertex $v$. We can describe $\BU(v)$ using the combinatorial framework for graph blow-ups described in \cite{culler86} and \cite{vogtmann90}, namely $\BU(v)$ is the poset of \emph{compatible partitions} of the set of incident half-edges, which we now recall.

\subsubsection*{Compatible partitions} Let $[m]\defeq \{1,2,\dots, m\}$, and consider partitions of $[m]$ into two blocks. Denote such a partition by $\alpha=\{a,\bar{a}\}$, where $1\in a$. Define the \emph{size} of $\alpha$ be
$$s(\alpha)\defeq |\bar{a}| \text{.}$$
Recall that distinct partitions $\{a,\bar{a}\}$ and $\{b,\bar{b}\}$ are said to be \emph{compatible} if either $a\subset b$ or $b\subset a$. For $m\geq3$ let $\Sigma(m)$ denote the simplicial complex of partitions $\alpha=\{a,\bar{a}\}$ of $[m]$ into blocks $a$ and $\bar{a}$ such that $a$ and $\bar{a}$ each have at least two elements, so $2\leq s(\alpha)\leq m-2$. That is, the vertices of $\Sigma(m)$ are such partitions, and a $j$-simplex is given by a collection of $j+1$ distinct, pairwise compatible partitions. Note that $\Sigma(3)=\emptyset$. Also define a similar complex $\Sigma'(m)$ for $m\geq2$, identical to $\Sigma(m)$ except that we allow partitions $\alpha=\{a,\bar{a}\}$ with $|\bar{a}|=1$. We do not allow $|a|=1$ though, so for example $\Sigma'(2)=\emptyset$.

For $v\neq p$ with $m\defeq \val(v)$, fix a labeling $1,\dots,m$ of the half-edges at $v$. Then the geometric realization of $\BU(v)$ is isomorphic to the barycentric subdivision of $\Sigma(m)$. In other words, a blow-up at $v$ is encoded by a chain of compatible partitions. A single partition describes an \emph{ideal edge}, i.e., an edge blow-up at a vertex, and the blocks $a$ and $\bar{a}$ indicate which half-edges attach to which endpoints of the new edge. See \cite{culler86} and \cite{vogtmann90} for more details.

\subsubsection*{Separating blow-ups} Thanks to Corollary~\ref{which_blowups_desc} we know precisely when a blow-up at $v\in\Lambda_i$ is descending, namely when it increases the number of vertices in $\Lambda_i$. Hence a blow-up at $v$ is descending if and only if it separates the set of half-edges at $v$ whose top is equal to $v$. We say that such a blow-up \emph{separates at} $v$. Let $\SBU(v)$ be the poset of blow-ups at $v$ that separate at $v$. Note that blow-ups at the basepoint $p$ are never separating, so $\SBU(p)=\emptyset$.

\subsubsection*{Splitting partitions} We will say that a partition $\alpha=\{a,\bar{a}\}$ of $[m]$ \emph{splits} a subset $S\subseteq[m]$ if $S\not\subseteq a$ and $a\not\subseteq S$. Define the \emph{splitting level} $\ell(\alpha)$ to be the minimum element of $\bar{a}$, i.e., the smallest $\ell$ such that $\alpha$ splits $[\ell]$. Note that $2\leq\ell(\alpha)\leq m-1$ for $\alpha\in\Sigma(m)$ and $2\leq\ell(\alpha)\leq m$ for $\alpha\in\Sigma'(m)$. Let $\Sigma(m,r)$ be the sublevel set of $\Sigma(m)$ spanned by partitions $\alpha$ with $\ell(\alpha)\leq r$, and similarly define $\Sigma'(m,r)$.

The next lemma gives a reformulation of $\Sigma(m,r)$ in terms of graph blow-ups. We assume now that in our fixed labeling of the half-edges of $v$, those half-edges whose top is $v$, say there are $r$ of them, are labeled precisely by $1,\dots,r$.

\begin{lemma}[Separating blow-ups and splitting partitions]\label{lem:blowups_and_partitions}
 Let $v\neq p$ be a vertex in $\Gamma$ with $m$ incident half-edges. Let $r$ be the number of half-edges with top $v$. Then $|\!\SBU(v)|\simeq \Sigma(m,r)$.
\end{lemma}

\begin{proof}
 The geometric realization $|\!\SBU(v)|$ contains the barycentric subdivision of $\Sigma(m,r)$ as a subcomplex. Also, any simplex in $|\!\SBU(v)|$ has at least one vertex in $\Sigma(m,r)$. Hence there is a map $|\!\SBU(v)|\to|\!\SBU(v)|$ sending each simplex to its face spanned by vertices in $\Sigma(m,r)$. This induces a deformation retraction from $|\!\SBU(v)|$ to $\Sigma(m,r)$.
\end{proof}

We now want to calculate the homotopy type of $\Sigma(m,r)$, and perhaps unsurprisingly we will use Morse theory. Consider the height function
$$z(\alpha)\defeq (\ell(\alpha),s(\alpha))$$
on $\Sigma(m)$, with the lexicographic ordering. Since compatible partitions have different sizes, they also have different $z$-values. Note that $\Sigma(m,r)$ is a sublevel set with respect to $z$, namely $\Sigma(m,r)=\Sigma(m)^{z\leq(r,m-2)}$. Hence we can analyze the homotopy type of $\Sigma(m,r)$ by looking at descending links in $\Sigma(m)$ with respect to $z$. We can also think of $z$ as a height function on $\Sigma'(m)$, and before handling $\Sigma(m,r)$ it will be convenient to first calculate the homotopy type of $\Sigma'(m,r)$.

\begin{lemma}
 For any $m\geq2$ and $2\leq r\leq m$, $\Sigma'(m,r)\simeq\bigvee S^{m-3}$.
\end{lemma}

\begin{proof}
 We induct on $m$. Since $\Sigma'(2)=\emptyset$, we already know that $\Sigma'(2,r)=\emptyset=S^{2-3}$ for any $r$, which handles the base case. Now let $m>2$ and consider the complex $\Sigma'(m,2)$. This is spanned by partitions $\{a,\bar{a}\}$ in which the set $\{1,2\}$ is split, and so any such $a$ will be $a=\{1\}\cup T$ for $T$ a non-empty subset of $\{3,4,\dots,m\}$. Thus $\Sigma'(m,2)$ is isomorphic to the barycentric subdivision of an $(m-3)$-simplex, and so is contractible.

 We now analyze the descending links of partitions with respect to $z$. Let $\alpha=\{a,\bar{a}\}$ be a partition in $\Sigma'(m,r)\setminus\Sigma'(m,2)$ and set $\ell\defeq \ell(\alpha)>2$ and $s\defeq s(\alpha)$. A partition $\beta=\{b,\bar{b}\}$ compatible with $\alpha$ is in the $z$-descending link $\dlk_z(\alpha)$ of $\alpha$ precisely when either $\ell(\beta)<\ell$, or $\ell(\beta)=\ell$ and $a\subsetneq b$. Note that in the first case $b\subseteq a$, so any partition of the first type is compatible with every partition of the second type. Hence the $z$-descending link of $\alpha$ is a join, of a \emph{d-in-link} and a \emph{d-out-link}. The d-in-link is the full subcomplex of $\dlk_z(\alpha)$ spanned by partitions of the first type, and the d-out-link is spanned by partitions of the second type. See Figure~\figref{in_and_out} for an example.

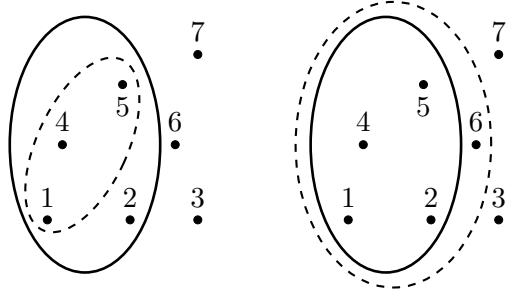
\begin{figure}[htb]\hypertarget{in_and_out}\leavevmode
    \centering
    
\begin{tikzpicture}
  \filldraw
  (0,0) circle (1.5pt)
  (1.1,0) circle (1.5pt)
  (2,0) circle (1.5pt)
  (0.2,1) circle (1.5pt)
  (1.7,1) circle (1.5pt)
  (1,1.8) circle (1.5pt)
  (2,2.2) circle (1.5pt);
  \node at (0,0.3) {$1$};
  \node at (1.1,0.3) {$2$};
  \node at (2,0.3) {$3$};
  \node at (0.2,1.3) {$4$};
  \node at (1.7,1.3) {$6$};
  \node at (1,1.5) {$5$};
  \node at (2,2.5) {$7$};
  \draw[line width=1pt] (0.5,1) ellipse (1 cm and 1.7 cm);
  \draw[line width=0.8pt, dashed, rotate=-25] (0,1.1) ellipse (0.6 cm and 1.25 cm);

  \begin{scope}[xshift=4cm]
   \filldraw
   (0,0) circle (1.5pt)
   (1.1,0) circle (1.5pt)
   (2,0) circle (1.5pt)
   (0.2,1) circle (1.5pt)
   (1.7,1) circle (1.5pt)
   (1,1.8) circle (1.5pt)
   (2,2.2) circle (1.5pt);
   \node at (0,0.3) {$1$};
   \node at (1.1,0.3) {$2$};
   \node at (2,0.3) {$3$};
   \node at (0.2,1.3) {$4$};
   \node at (1.7,1.3) {$6$};
   \node at (1,1.5) {$5$};
   \node at (2,2.5) {$7$};
   \draw[line width=1pt] (0.5,1) ellipse (1 cm and 1.7 cm);
   \draw[line width=0.8pt, dashed] (0.6,1) ellipse (1.3 cm and 1.9 cm);
  \end{scope}

\end{tikzpicture}

    \caption{A partition in the d-in-link, and one in the d-out-link, of a partition with size $s=3$ and splitting level $\ell=3$.}
    \label{fig:in_and_out}
\end{figure}

 First consider the d-out-link. Partitions $\beta=\{b,\bar{b}\}$ in the d-out-link are characterized by the property that $a\subsetneq b$ and $\ell\in\bar{b}$. Treating $a$ as a single point, this amounts to saying that $a\subsetneq b$ and $\beta$ splits $\{a,\ell\}$. Hence the d-out-link is isomorphic to $\Sigma'(s+1,2)$. If $s=1$ this is empty, and if $s>1$ this is contractible as explained above. In particular if $s>1$ then $\dlk_z(\alpha)$ is already contractible. Now assume $s=1$, so the d-out-link is empty and $\dlk_z(\alpha)$ just equals the d-in-link. Then the d-in-link is isomorphic to the complex of partitions of $[m-1]$ that split $[\ell-1]$, and so is given by $\Sigma'(m-1,\ell-1)$. This is $(m-1-3)$-spherical by induction, so we conclude that all descending links are either contractible or $(m-4)$-spherical. Since $\Sigma'(m,2)$ is $(m-3)$-spherical this implies that $\Sigma'(m,r)$ is also $(m-3)$-spherical \cite[Corollary~2.6]{bestvina97}.
\end{proof}

\begin{proposition}\label{prop:desc_partitions_conn}
 For any $m\geq3$ and $2\leq r\leq m-1$, $\Sigma(m,r)\simeq\bigvee S^{m-4}$.
\end{proposition}

\begin{proof}
 As in the previous proof we induct on $m$. When $m=3$ we only consider $r=2$, and $\Sigma(3,2)$ is empty. Now let $m>3$ and consider $\Sigma(m,2)$. As with $\Sigma'(m,2)$, $\Sigma(m,2)$ is spanned by partitions $\{a,\bar{a}\}$ in which the set $\{1,2\}$ is split, and so any such $a$ will be $a=\{1\}\cup T$, for $T$ now a \emph{proper} non-empty subset of $\{3,4,\dots,m\}$. (Now we cannot have $T=\{3,4,\dots,m\}$ since the resulting partition would have size $1$.) Thus $\Sigma(m,2)$ is the surface of a barycentrically subdivided $(m-3)$-simplex, and so is homeomorphic to $S^{m-4}$.
 
 Now consider the descending link $\dlk_z(\alpha)$ of $\alpha=\{a,\bar{a}\}$ with $\ell\defeq \ell(\alpha)>2$ and $s\defeq s(\alpha)$. The descending link decomposes as before as the join of a d-in-link and d-out-link. By the same argument as in the previous proof, the d-out-link is isomorphic to $\Sigma(s+1,2)$, which is homeomorphic to $S^{s-3}$. The d-in-link is isomorphic to the complex of partitions of $[m-s]$ that split $[\ell-1]$ and have size at least $1$. (Since $\bar{a}$ has elements in it, we do have to consider partitions of $[m-s]$ that have size $1$ as a partition of $[m-s]$.) So, the d-in-link is isomorphic to $\Sigma'(m-s,\ell-1)$, and hence is homotopy equivalent to $\bigvee S^{m-s-3}$ by the previous lemma. Then $\dlk_z(\alpha)$ is the join of the d-in- and d-out-links, and so is homotopy equivalent to $(\bigvee S^{m-s-3})\ast S^{s-3}=\bigvee S^{m-5}$. Since $\Sigma(m,2)$ is $(m-4)$-spherical and the descending links of partitions in $\Sigma(m,r)\setminus\Sigma(m,2)$ are all $(m-5)$-spherical, we conclude that $\Sigma(m,r)$ is $(m-4)$-spherical \cite[Corollary~2.6]{bestvina97}.
\end{proof}

We remark that since $\Sigma(m,m-1)=\Sigma(m)$, we recover the fact that $\Sigma(m)$ is $(m-4)$-spherical, as shown in \cite[Theorem~2.4]{vogtmann90}. Coupling Proposition~\ref{prop:desc_partitions_conn} with Lemma~\ref{lem:blowups_and_partitions} we see that if there are least two half-edges with top $v$, then
$$|\!\SBU(v)|\simeq\bigvee S^{\val(v)-4} \text{.}$$

Now let $\displaystyle A\defeq \ast_{v\neq p}\SBU(v)$, where the join is taken over all vertices $v\neq p$ in $\Gamma$. Recall that $V$ is the number of vertices in $\Gamma$.

\begin{corollary}\label{cor:upcor}
 If $\Gamma$ has no decisive edges then $|A|\simeq\bigvee S^{d_0(\Gamma)-V}$.
\end{corollary}

\begin{proof}
 Since there are no decisive edges, for any $v\neq p$ we know that there are at least two half-edges at $v$ with top $v$. Hence $|\!\SBU(v)|\simeq\bigvee S^{\val(v)-4}$, and so
\[|A|\simeq\ast_{v\neq p}\bigvee S^{(\val(v)-2)-2}=\bigvee S^{(d_0(\Gamma)-2(V-1))+(V-2)}=\bigvee S^{d_0(\Gamma)-V} \text{.}\qedhere\]
\end{proof}

\begin{proposition}[Homotopy type of the d-up-link]\label{prop:upprop}
 If $\Gamma$ has no decisive edges then the d-up-link is homotopy equivalent to $|A|$, and hence to $\bigvee S^{d_0(\Gamma)-V}$.
\end{proposition}

\begin{proof}
 For a poset $P$, define $\underline{P}$ to be $P\sqcup\{\bot\}$, with $\bot$ a formal minimum element. Then $P\ast Q\cong\underline{P}\times\underline{Q}\setminus\{(\bot,\bot)\}$ for posets $P$ and $Q$. The relevant example is that
 $$A=\ast_{v\neq p}\SBU(v)\cong\prod_{v\neq p}\underline{\SBU}(v)-\{(\bot)_v\}\eqdef Y \text{.}$$
Define \[\displaystyle X\defeq \biggl\{f\in\prod_{v\neq p}\underline{\BU}(v)\biggm\vert\exists v\in\Lambda_{D(f)}\text{ with }f_v\in\SBU(v)\biggr\}.\] Here $f_v$ is the blow-up at vertex $v$ in the tuple $f$, and $D(f)$ is the minimal level such that $f_v\neq\bot$ for some $v\in\Lambda_{D(f)}$. Note that $Y\subseteq X$. Define a map $r \colon X\to X$ by

\begin{center}
$(f_v)_v\mapsto\left(\begin{cases}
  f_v & \text{for} \ f_v \in \SBU(v)\\
  \bot & \text{for} \ f_v \not\in \SBU(v)
\end{cases}
\right)_v$.
\end{center}

 Note that $r$ is a poset map that is the identity on its image $Y$. Also, $r(f)\leq f$ for all $f\in X$, so $r$ induces a homotopy equivalence between $|X|$ and $|Y|$ \cite[Section~1.3]{quillen78}. But $|X|$ is precisely the d-up-link of $\Gamma$, so the d-up-link is homotopy equivalent to $\bigvee S^{d_0(\Gamma)-V}$ by Corollary~\ref{cor:upcor}.
\end{proof}

\section{Proof of the main results}\label{sec:main_proofs}

\begin{corollary}[Homotopy type of descending links]\label{cor:desc_lk_conn}
 For any vertex $\Gamma$ in $L_n$, $\dlk(\Gamma)$ is either contractible or homotopy equivalent to $\bigvee S^{d_0(\Gamma)-1}$.
\end{corollary}

\begin{proof}
 If the d-down-link of $\Gamma$ is contractible, then so is $\dlk(\Gamma)$. If the d-down-link is not contractible, then $\Gamma$ has no decisive edges (Lemma~\ref{lem:downlemma2}). Hence joining the d-up-link and d-down-link yields \[\displaystyle\left(\bigvee S^{d_0(\Gamma)-V}\right)\ast\left(\bigvee S^{V-2}\right)\simeq\bigvee S^{d_0(\Gamma)-1}\] (Propositions~\ref{prop:downprop} and~\ref{prop:upprop}).
\end{proof}

\begin{theorem}[Degree Theorem]\label{thrm:deg_thm}
 $L_{n,k}$ is $(k-1)$-connected.
\end{theorem}

\begin{proof}
 For any vertex $\Gamma$ in $L_n\setminus L_{n,k}$ we have $d_0(\Gamma)>k$, so by the previous corollary, $\dlk(\Gamma)$ is $(k-1)$-connected. Since $L_n$ is contractible and $L_{n,k}$ is a sublevel set of $L_n$ with respect to $h$ (Observation~\ref{obs:sublev}), $L_{n,k}$ is $(k-1)$-connected by \cite[Corollary~2.6]{bestvina97}.
\end{proof}

\end{document}